\documentclass{amsart}

\usepackage{amsmath,amssymb,amsthm}

\usepackage{tikz}
\usetikzlibrary{calc,intersections,through,backgrounds}
\usepackage[latin1]{inputenc}

\usepackage{caption}
\usepackage{graphicx}
\usepackage{graphics}

\newtheorem{theorem}{Theorem}

\newtheorem*{lemma}{Lemma}

\begin{document}

\title[]{Regularized Potentials of Schr{\"o}dinger Operators\\ and a local landscape function}
\author[]{Stefan Steinerberger}
\address{Stefan Steinerberger, Department of Mathematics, Yale University, 10 Hillhouse Avenue, New Haven, CT 06511, USA} \email{stefan.steinerberger@yale.edu}
\thanks{The author is by the NSF (DMS-1763179) and the Alfred P. Sloan Foundation.}

\begin{abstract} We study localization properties of low-lying eigenfunctions $$(-\Delta +V) \phi = \lambda \phi \qquad \mbox{in}~\Omega$$ for rapidly varying potentials $V$ in bounded domains $\Omega \subset \mathbb{R}^d$. 
 Filoche \& Mayboroda introduced the landscape function $(-\Delta + V)u=1$ and showed that the function $u$ has remarkable properties: localized eigenfunctions
 prefer to localize in the local maxima of $u$. Arnold, David, Filoche, Jerison \& Mayboroda showed that $1/u$ arises naturally as the potential in a related equation.
Motivated by these questions, we introduce
a one-parameter family of regularized potentials $V_t$ that arise from convolving $V$ with the radial kernel 
$$ V_t(x) = V * \left(  \frac{1}{t} \int_0^t \frac{ \exp\left( - \|\cdot\|^2/ (4s) \right)}{(4 \pi s )^{d/2}} ds \right).$$ 
We prove that for eigenfunctions $(-\Delta +V) \phi = \lambda \phi$ this regularization $V_t$ is, in a precise sense, the canonical effective potential on small scales. The landscape function $u$ respects the same type of regularization. This allows allows us to derive landscape-type functions out of solutions of the equation $(-\Delta + V)u = f$ for a general right-hand side $f:\Omega \rightarrow \mathbb{R}_{>0}$.
\end{abstract}
\maketitle

\vspace{-10pt}

\section{Introduction}
\subsection{The Landscape function.} Physical systems comprised of inhomogeneous materials sometimes exhibit localized vibration patterns: throughout this paper, let $\Omega \in \mathbb{R}^d$ be an open, bounded domain in which we consider the equation
\begin{align*}
(-\Delta + V)\phi &= \lambda \phi~ \quad \mbox{in~}\Omega \\
 \phi&= 0 \qquad \mbox{on}~ \partial \Omega
\end{align*}
where $V:\Omega \rightarrow \mathbb{R}_{\geq 0}$ is a real-valued, nonnegative potential. If $V$ oscillates rapidly, then this equation may have eigenfunctions
that are strongly localized \cite{anderson}. These determine the behavior of many associated dynamical systems (say, the heat equation, the wave equation or the Schr{\"o}dinger equation) and
are of obvious interest.
Filoche \& Mayboroda \cite{fil} have provided a simple but astonishingly effective method to predict the behavior of low-energy eigenfunctions for such operators $-\Delta + V$. They define the \textit{landscape function} as the unique function $u:\Omega: \mathbb{R} \rightarrow \mathbb{R}_{}$ solving
\begin{align*}
(-\Delta + V)u &=1~ \qquad \mbox{in~}\Omega \\
 u&= 0 \qquad  \mbox{on}~~ \partial \Omega
\end{align*}
and show that $u$ exerts pointwise control on all eigenfunctions $(-\Delta + V)\phi = \lambda \phi$ 
$$ |\phi(x)|  \leq \lambda u(x)  \|\phi\|_{L^{\infty}(\Omega)} .$$ 
An eigenfunction $\phi$ can only localize in $\left\{x: u(x) \geq 1/\lambda\right\} \subset \Omega$. However, the landscape function turns out to be much more effective than that. Numerical experiments \cite{fil} suggest that the largest local maxima correspond precisely to the location where the first few eigenfunctions localize and that many more properties (including refined eigenvalue estimates and improvements on the Weyl law) are being captured. The accuracy of these refined predictions is quite striking and have already led to many interesting results \cite{arnold0, arnold, arnold2, chal, david, fil, fil2, fil3, har, lef, leite, pic}.

\begin{center}
\begin{figure}[h!]
\begin{tikzpicture}
\node at (0,0) {\includegraphics[width=0.6\textwidth]{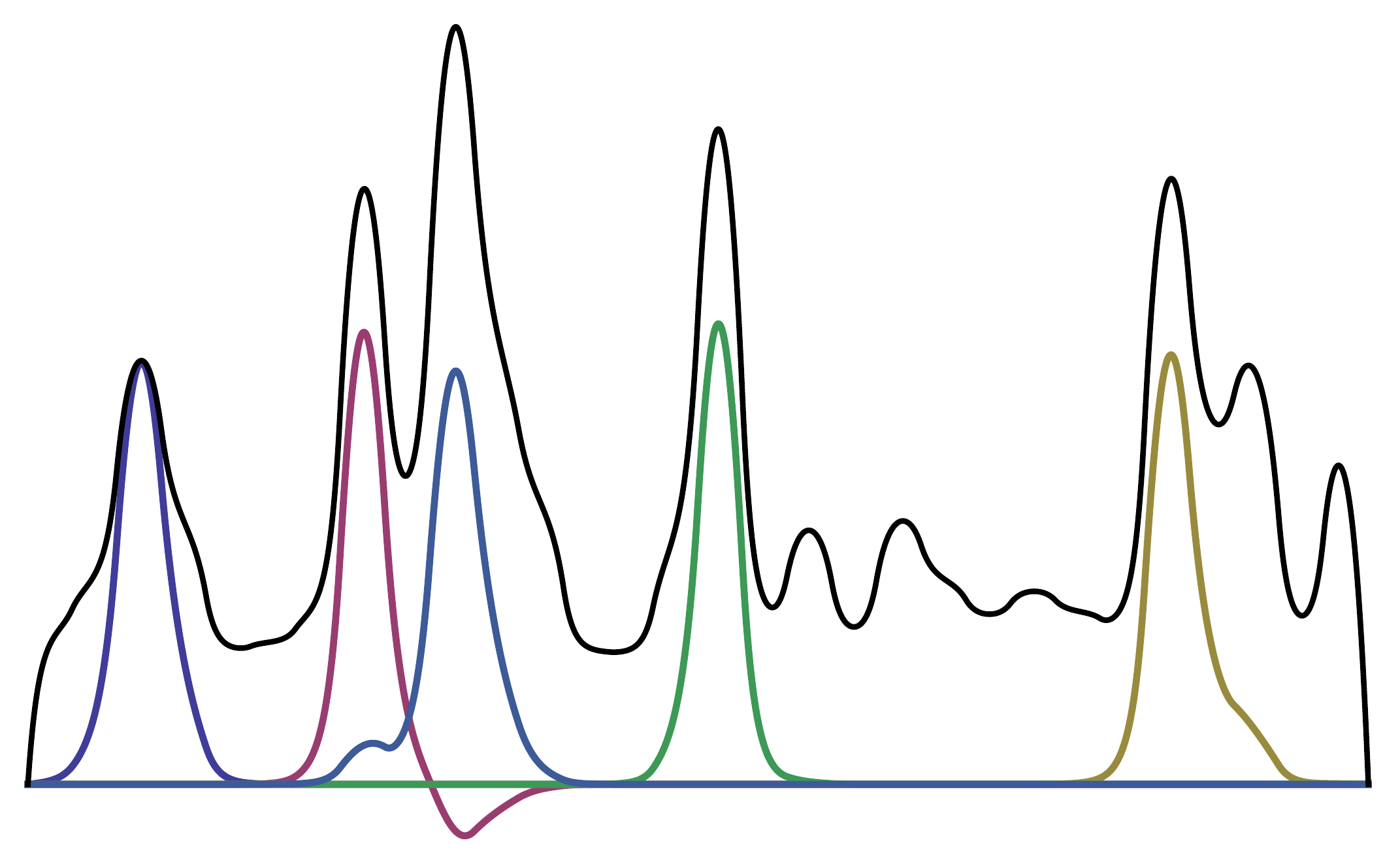}};
\end{tikzpicture}
\caption{An example of $(-\Delta + V)\phi = \lambda \phi$ on the unit interval with a rapidly oscillating and irregular potential $V$. The peaks in the Filoche-Mayboroda landscape function (black) predict where the first five eigenfunctions (in color) localize.}
\end{figure}
\end{center}

\vspace{-20pt}

\subsection{The effective potential.} 
$1/u$ seems to be better than $V$ when it comes to describing the localization properties of the eigenfunctions \cite{arnold}. There is a simple reason why one might expect `effective potentials', potentials derived from $V$, to have more predictive power than $V$ itself: the eigenfunctions of $-\Delta + V$ try to minimize the $L^2-$norm of their gradient,
$\|\nabla \phi\|_{L^2}^2$, while simultaneously trying to minimize $\left\langle \phi, V \phi \right\rangle$. This competition leads to the eigenfunction $\phi$ not really `seeing' $V$ but rather `seeing'
an averaged, smoothed or somehow regularlized version of $V$.
One such an effective potential was proposed by
 Arnold, David, Jerison, Mayboroda \& Filoche \cite{arnold}. Their approach is based on writing an eigenfunction as $\phi = u \psi$ for some unknown function $\psi$. The equation
$$ \left( - \Delta + V \right) \phi = \lambda \phi$$
then transforms into
$$ - \frac{1}{u^2} \mbox{div}(u^2 \nabla \psi)  +  \frac{1}{u} \psi = \lambda \psi.$$
The new dominating potential $1/u$ is now responsible for the underlying dynamics for this related equation. In particular, this reformulation allows for Agmon estimates. Introducing an Agmon distance
$$ \rho(r_1, r_2) = \min_{\gamma}\left( \int_{\gamma}{\sqrt{(\omega(r) - \lambda)_{+}} ds} \right),$$
where $\gamma$ ranges over all paths from $r_1$ to $r_2$
and using Agmon's inequality \cite{agmon}, one can deduce that for eigenfunctions $\phi$ localized in $r_0 \in \Omega$ 
$$ |\phi(r)| \lesssim e^{-\rho(r_0, r)}.$$
 There is convincing numerical evidence \cite{arnold} that $1/u - \lambda$
seems to predict decay more accurately than the classical quantity $V-\lambda$. This might seem surprising because $V$ determines the behavior of the eigenfunctions by being the term arising in the equation and is again due to the phenomenon described above: eigenfunctions are forced to be close to constant on small scales which leads to them interacting with a regularized version of $V$ instead of $V$ itself.

\subsection{Related results.} Other methods for the purpose of fast computation of the location of localized low-lying eigenfunctions have been proposed \cite{alt, altmann, korn, lu}. A first attempt at a local description of the landscape function was given by the author in \cite{steini}.  We also mention a curious localization phenomenon for Neumann boundary conditions \cite{sapo1, jones,  sapo2}.

\section{The Result}
\subsection{Introduction.} This section presents the main idea. We summarize the existing insights which motivate our approach.
\begin{enumerate}
\item If the potential $V$ is smooth and slowly varying (say, essentially constant on scales larger than the wavelength $\lambda^{-1/2}$), then low-lying eigenfunctions 
$$ (-\Delta + V)\phi = \lambda \phi$$
localize in the local minima of $V$ (this is, in a certain sense, the regime of classical physics).
\item  If the potential $V$ is rough, irregular and quickly varying (say, oscillating dramatically on scales comparable to or smaller than the wavelength $\lambda^{-1/2}$), then the requirement of keeping $\|\nabla \phi\|_{L^2}^2$ small starts playing a more fundamental role: localized eigenfunction do not interact with $V$ as much as they interact with a locally regularized version of $V$ (this is, in a certain sense, the regime of quantum physics).
\item The remaining question is: \textit{what} is the regularized potential? Arnold, David, Jerison, Mayboroda \& Filoche \cite{arnold} show that $1/u$ is a possible regularization of $V$ (arising as the potential for a related equation).\\
\end{enumerate}

One natural question is whether localized eigenfunctions $(-\Delta + V)\phi = \lambda \phi$, while seemingly not directly interacting with $V$, perhaps interact with a locally averaged version of $V$. This seems natural when considering the eigenfunction $\phi$ as a critical point of the Dirichlet energy (subject to orthogonality to previous eigenfunctions)
$$J(\phi) =  \int_{\Omega}{ |\nabla \phi|^2 dx} + \int_{\Omega}{ V \phi^2 dx}.$$ 
 The first term is of a certain size: in particular, this forces eigenfunctions to be essentially constant below a certain scale. However, if they are constant over small scales, then the second integral really assumes a very different meaning: what is relevant is not the value of $V$ so much as the local average of $V$ (averaged over the scale over which we expect $\phi$ to be constant).
The question then naturally is: what would be a natural way of averaging $V$? 
The main point of our paper is to show that there is a \textit{canonical} way of computing local averages of $V$ in a way that respects the behavior of eigenfunctions in a precise sense.

\subsection{A Convolution Kernel.} We define the local average $V_t$ of the potential $V$ at a scale $t>0$ as the convolution of the potential $V$ with the
radially symmetric kernel $k_t:\mathbb{R}^d \rightarrow \mathbb{R}_{\geq 0}$ that is given by
$$ k_t(x) = \frac{1}{t} \int_0^t \frac{ \exp\left( - \|x\|^2/ (4s) \right)}{(4 \pi s )^{d/2}} ds.$$
The kernel depends on a scale parameter $t>0$ and
the dimension $d$ of the domain $\Omega \subset \mathbb{R}^d$ but nothing else. The radial profiles of the kernel in $d=1$ and $d=2$ dimensions are shown in Fig. 2. 
\begin{center}
\begin{figure}[h!]
\begin{tikzpicture}
\node at (0,0.1) {\includegraphics[width=0.42\textwidth]{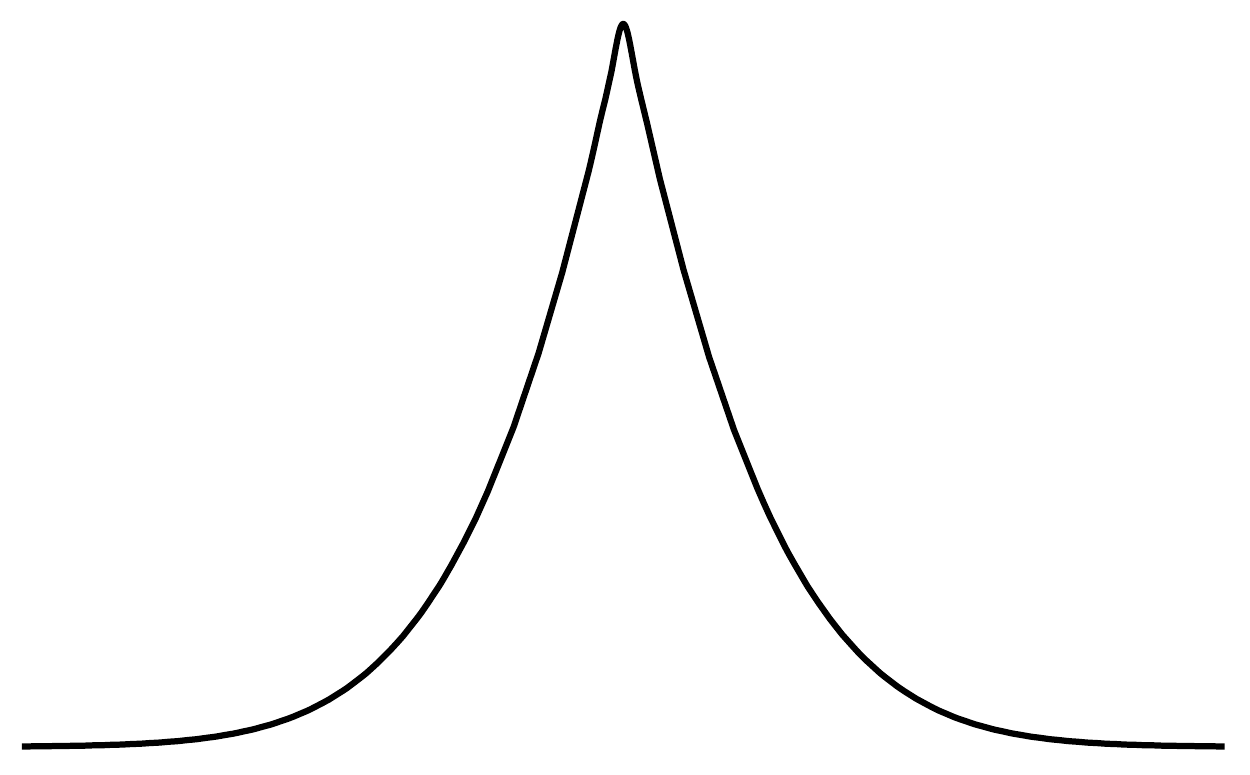}};
\draw [<->] (-3,-1.52) -- (3,-1.52);
\draw [->] (-0.01,-1.52) -- (-00.01,2);
\node at (0., -1.7) {0};
\node at (7,0.1) {\includegraphics[width=0.42\textwidth]{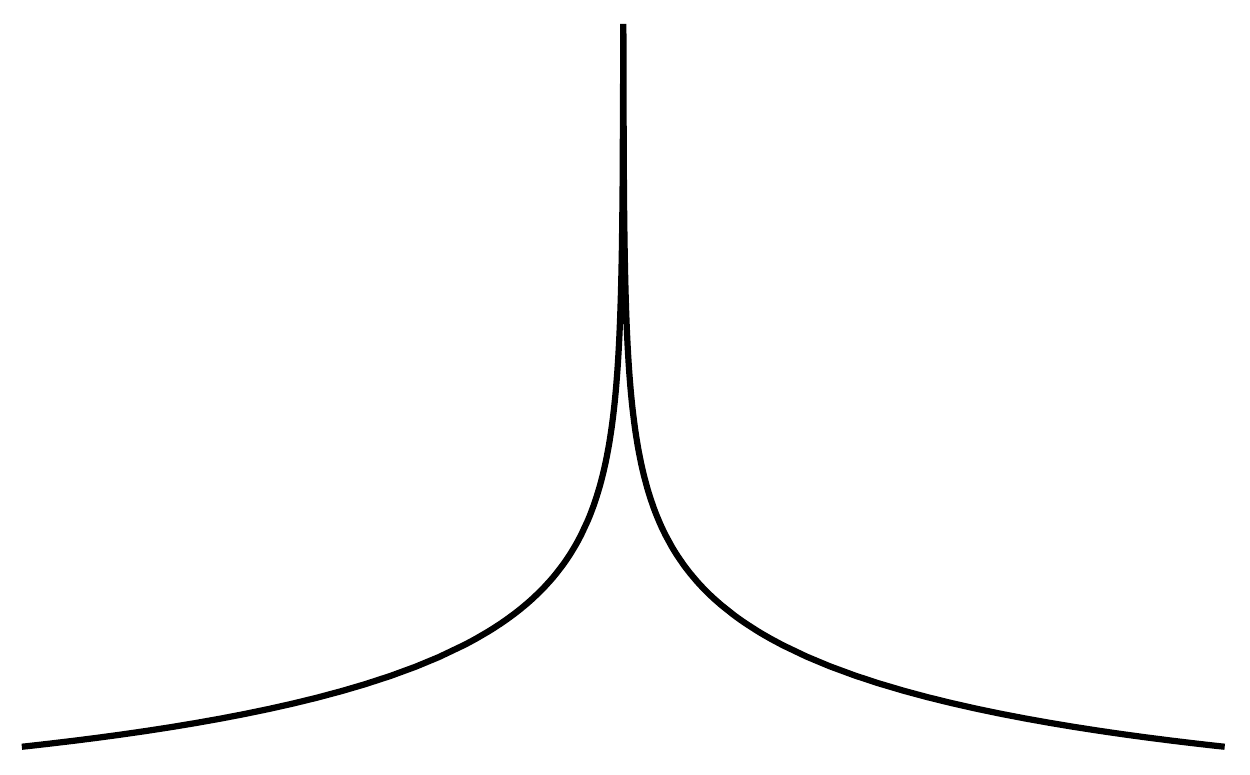}};
\draw [<->] (7-3,-1.52) -- (7+3,-1.52);
\draw [->] (6.99,-1.52) -- (6.99,2);
\node at (7, -1.7) {0};
\end{tikzpicture}
\caption{The radial profile of the convolution kernel $k_t(r)$ in $d=1$ dimensions (left) and $d=2$ (right) that we derive below.}
\end{figure}
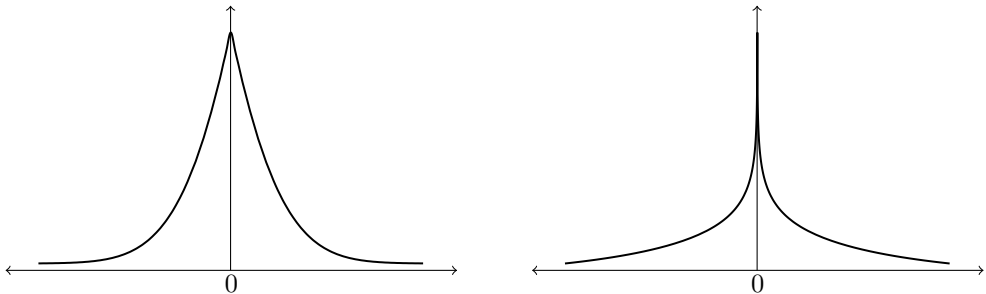
\end{center}
These kernels have different closed forms in different dimensions, for example
\begin{align*}
k_t(r) &= \frac{1}{\sqrt{\pi t}} \exp\left( -\frac{r^2}{4t}\right) - \frac{r}{2t} \mbox{erfc}\left(\frac{r}{2\sqrt{t}}\right) \qquad \qquad (d=1)\\
k_t(r) &= \frac{1}{4\pi t}\Gamma\left(0, \frac{r^2}{4t}\right)  \qquad \quad \qquad \qquad \qquad \qquad \qquad (d=2),
\end{align*}
where $\Gamma(0,z)$ is the incomplete gamma function. We observe that $k_t$ has normalized total mass, i.e. in all dimensions and for all $t>0$,
$$ \int_{\mathbb{R}^n} \left( \frac{1}{t} \int_0^t \frac{ \exp\left( - \|x\|^2/ (4s) \right)}{(4 \pi s )^{d/2}} ds \right) dx = 1.$$
Most of the $L^1-$mass of $k_t$ is concentrated at scale $\sim \sqrt{t}$ around the origin and $k_t$ is exponentially decaying after that, this follows easily from observing that it is a linear combination of Gaussians the widest of which is $\exp(-\|x\|^2/(4t)$. These kernels $k_t$ approximate the identity as $t \rightarrow 0$ in the
sense that $(V * k_t)(x) \rightarrow V(x)$ whenever $V$ is bounded and continuous in a neighborhood of $x$. \\

 In particular, which is shown in simple numerical examples throughout the paper, the convolution $V_t = V * k_t$ does appear
to behave like a one would expect from a regularized potential: low-lying eigenfunctions minimize in the local minima of $V_t$ (whereas the local minima of $V$ have relatively little explanatory power). It seems to track $1/u$ quite closely (see also Theorem 2 and \S 2.6).

\subsection{The Result} We suppose we are given an eigenfunction
\begin{align*}
(-\Delta + V)\phi &= \lambda \phi~ \quad \mbox{in~}\Omega \\
 \phi&= 0 \qquad \mbox{on}~ \partial \Omega.
\end{align*}
We will now try to understand how the solution of this equation behaves under convolving $V$ with a kernel. We want to average $V$ over as large as possible a region while
still almost satisfying the equation; this naturally identifies the kernel $k_t$ (which, indeed, is derived from the proof of the Theorem). The Theorem is somewhat subtle:
the crux of the statement is that an implicit constant only depends on the size of the potential $V$ but not on any of its derivatives or finer properties.

\begin{center}
\begin{figure}[h!]
\begin{tikzpicture}
\node at (0,0) {\includegraphics[width=0.7\textwidth]{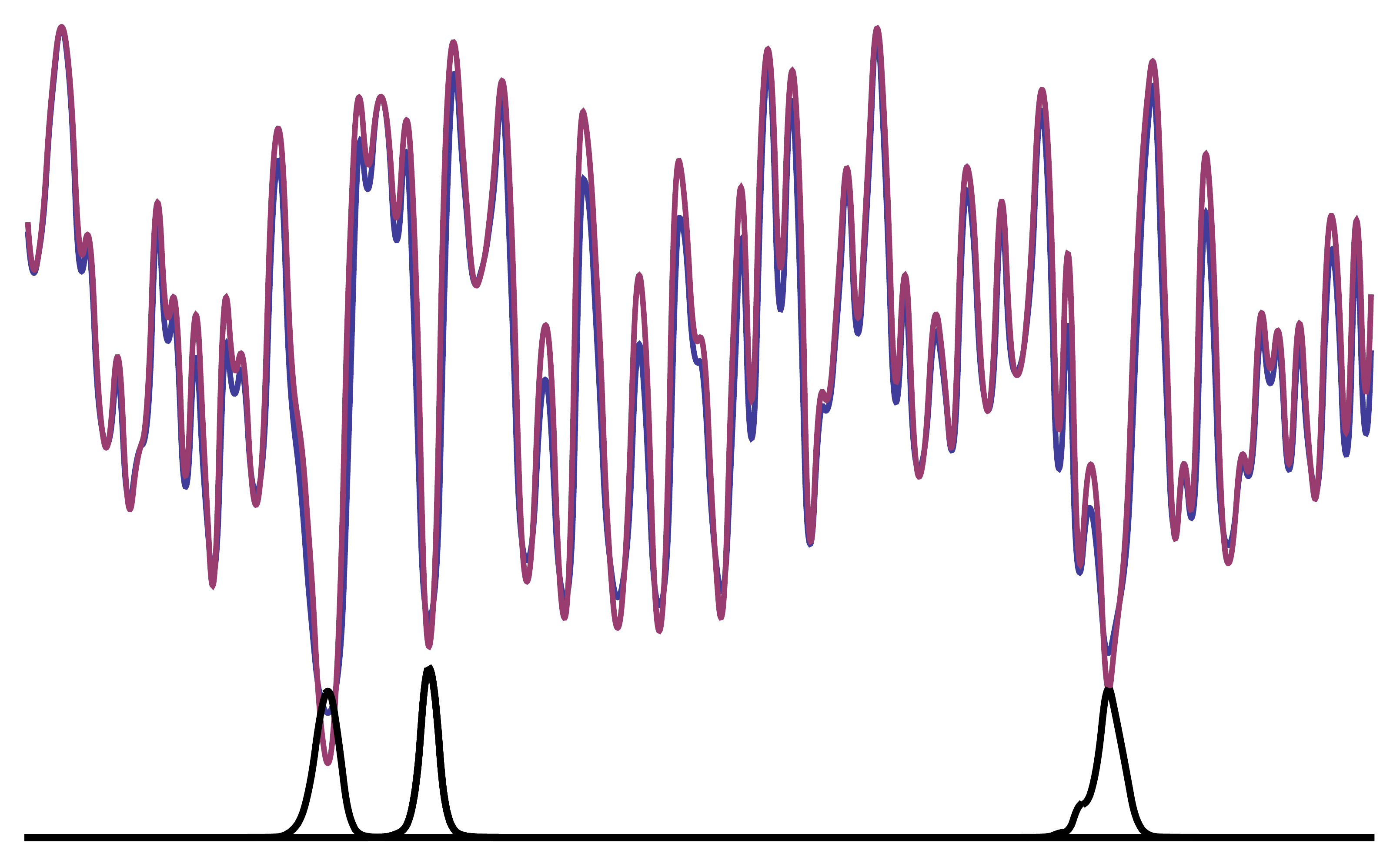}};
\end{tikzpicture}
\caption{The inverse landscape function $1/u$ (blue) and the regularized potential $V_t = V* k_t$ (purple). The first three eigenfunctions (black) localize in the regularized potential wells.}
\end{figure}
\end{center}

\begin{theorem}[Main Result]
Let $\Omega \subset \mathbb{R}^n$ be an open, bounded domain with smooth boundary, let $0 \leq V \in C(\overline{\Omega})$ be a continuous potential and let $\phi$ be a solution of
\begin{align*}
(-\Delta + V)\phi &= \lambda \phi~ \quad \mbox{in~}\Omega \\
 \phi&= 0 \qquad \mbox{on}~ \partial \Omega.
\end{align*}
Then, for any fixed $x \in \Omega$, as $t \rightarrow 0$, we have, for $k_t$ as above,
$$ -\Delta \phi(x) +  (V*k_t)(x) \phi(x)= \lambda \phi(x) + \mathcal{O}_{\phi, \|V\|_{L^{\infty}}}(t),$$
where the implicit constant depends \emph{only} on $\phi$ and $\|V\|_{L^{\infty}}$. 
\end{theorem} 
The essence of this result is that the error term only depends on the size of the potential and none of the finer properties (something that
is not true for other convolution kernels), see below for more remarks on this. This means that the equation barely changes when replacing $V$ by $V_t$ and this motivates our interpretation of
$V_t$ being the first-order approximation of the kernel \textit{as seen by the eigenfunction}: even though $\phi$ is created by $V$, it behaves \textit{almost} as if it were created by $V_t$.\\

\textbf{Remarks.}\\

(1)  The distinguished role of the kernel $k_t$ is somewhat subtle: it has the special property that the implicit constant
in the error term does not depend on any fine properties of the potential $V$ but merely on its size. That this is a nontrivial property even for smooth potentials $V$ can be seen as follows: let $g_t$ be another
radial kernel at scale $\sim \sqrt{t}$, for example a suitably scaled Gaussian. Then
$$ (V*k_t)(x) - (V*g_t)(x) = V*(k_t - g_t) = \int_{\mathbb{R}^d} V(x+y)\left( k_t(y) - g_t(y) \right) dy.$$
The function $k_t(y) - g_t(y)$ has mean value 0 and is localized at scale $\sim \sqrt{t}$. A Taylor expansion of $V$ in $x$ shows that, since both are radial, for some
universal constant $c$ depending only on the kernels,
$$  \int_{\mathbb{R}^d} V(x+y)\left( k_t(y) - g_t(y) \right) dy \sim c \cdot \Delta V(x) t \qquad \mbox{as}~t \rightarrow 0.$$
$V$ need not even be differentiable. Considering a potential behaving like $V(x) = \|x-x_0\|^{\alpha}$ in $x_0$ shows
that it is possible for the error term to be size $\mathcal{O}(t^{\alpha/2})$ for any $\alpha>0$ since the kernels are localized at scale $\sim \sqrt{t}$.\\

(2) Our assumption of $V$ being continuous is for simplicity of exposition. It does not have an impact on the applicability of the result in practice: low-lying eigenfunctions do not distinguish between discontinuous
potentials and continuous potentials as long as the continuous potential is allowed to have very large derivatives. (In some sense, this is the main idea behind all these investigations of the landscape function in the first place: to find the simpler `effective' potential). Any restrictions on the size of $\nabla V$ or $\Delta V$ would
severely affect the applicability of the result. The whole point of Theorem 1 is that the regularization $V_t$ does not induce errors depending on any such quantities. Our argument
is somewhat flexible and would apply to discontinuous potentials $V$ as well: $V$ being in the Kato class is the natural limitation of our method, we comment on this after the proof.

\begin{center}
\begin{figure}[h!]
\begin{tikzpicture}
\node at (0,0) {\includegraphics[width=0.6\textwidth]{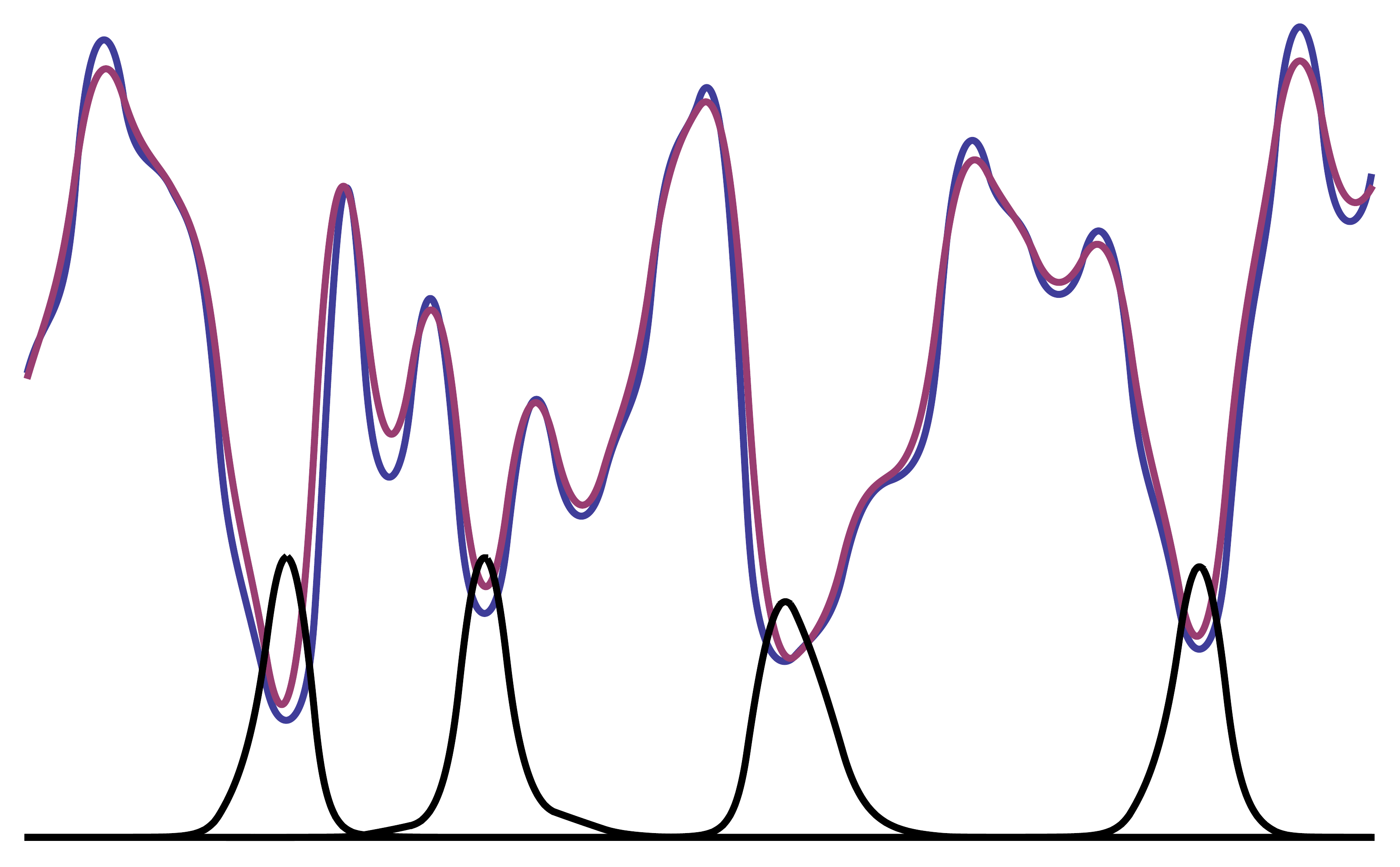}};
\end{tikzpicture}
\caption{ The inverse of the landscape function $1/u$ (blue) and the localized landscape $V_t = V* k_t$ (purple) are similar. The first eigenfunctions (black) localize in the regularized potential wells.}
\end{figure}
\end{center}

(3) The size of the error term can indeed be specified in terms of $\|V\|_{L^{\infty}}$, however, the proof actually implies a slightly refined statement. The dependence on $\|V\|_{L^{\infty}}$ could be replaced by
the size of $V$ averaged over the scale $\sim \sqrt{t}$ (which, in turn, is trivially dominated by $\|V\|_{L^{\infty}}$). We believe that this could also be a natural choice for $t$ in applications. In the same manner, the error term does not depend on $\phi$ in a strong global sense
but merely on the size of the second order derivatives in a neighborhood of $x$.

\begin{center}
\begin{figure}[h!]
\begin{tikzpicture}
\node at (0,0) {\includegraphics[width=0.6\textwidth]{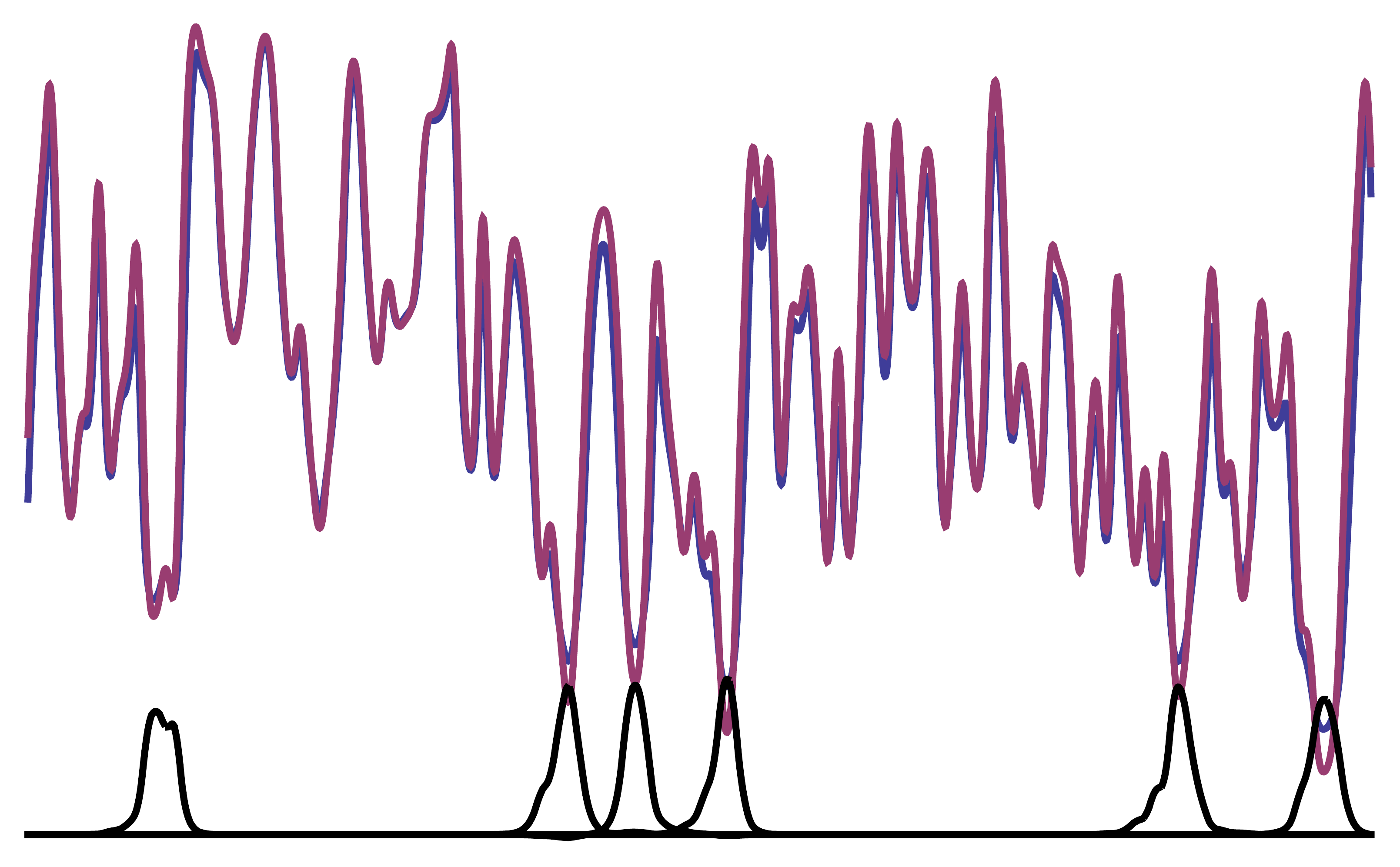}};
\end{tikzpicture}
\caption{The inverse landscape function $1/u$ (blue) and the localized landscape $V_t = V* k_t$ (purple). The first six eigenfunctions (black) localize in the regularized potential wells.}
\end{figure}
\end{center}

\subsection{Regularizing the Landscape Function.}
Throughout this paper, we consider simple numerical examples of the type
$$ (-\Delta + V)\phi = \lambda \phi \qquad \mbox{on}~[0,1]$$
with Dirichlet boundary conditions. We discretize the problem on $n=3000$ nodes. Our potentials $V$ are chose to be piecewise constants on short intervals (of size $\sim 1/20$ or $\sim 1/30$) and whose values are given by independently and identically distributed random values chosen uniformly at random from  $\sim [0, 10^5]$ (both parameters are slightly varied to produce different examples throughout the paper). As we can see in many examples, the inverse of the landscape function $1/u$ and the regularized potential $V_t$ behave remarkably similar in most cases.
This turns out to not be a coincidence. In fact, the landscape function exhibits a similar degree of stability under the type of regularization that we introduce in this paper.
\begin{theorem} 
Let $\Omega \subset \mathbb{R}^n$ be an open, bounded domain with smooth boundary, let $0 \leq V \in C^{}(\overline{\Omega})$ be a continuous potential. Let $x \in \Omega$, then
$$ -\Delta u(x) + (V*k_t)(x) u(x)  = 1 + \mathcal{O}_{u, \|V\|_{L^{\infty}}}(t),$$
where the error term depends \emph{only} on $u$ and $\|V\|_{L^{\infty}}$. 
\end{theorem}
As in Theorem 1, the remarkable aspect is that the error term does not depend on fine properties of $V$ around $x$ but merely on $\|V\|_{L^{\infty}}$ and the derivatives of $u$ itself.
Summarizing, we see that both the eigenfunctions of $-\Delta + V$
as well as the landscape function $u$ exhibit a strong form of stability under replacing $V$ by $V_t$. However, $V_t$ is a much nicer potential: one would expect that eigenfunctions localize
in local minima of $V_t$, say in $V_t(x_{\min})$ and that, because of the increased regularity, the solution of $(-\Delta + V_{t})u = 1$ satisfies $u(x_{\min}) \sim 1/V_t(x_{\min})$ because the
Laplacian plays a less dominant role. By the same reasoning, one would then expect that $1/u(x_{\min}) \sim \lambda \sim V_t(x_{\min})$.
 We believe this could conceivably play a role in the effectiveness of $1/u$, see also \S 2.6.
The connection between the landscape function and the
short-time asymptotic of an associated diffusion process, a tool crucial in our proof of Theorem 2, was already noted in \cite{steini}.

\begin{center}
\begin{figure}[h!]
\begin{tikzpicture}
\node at (0,0) {\includegraphics[width=0.6\textwidth]{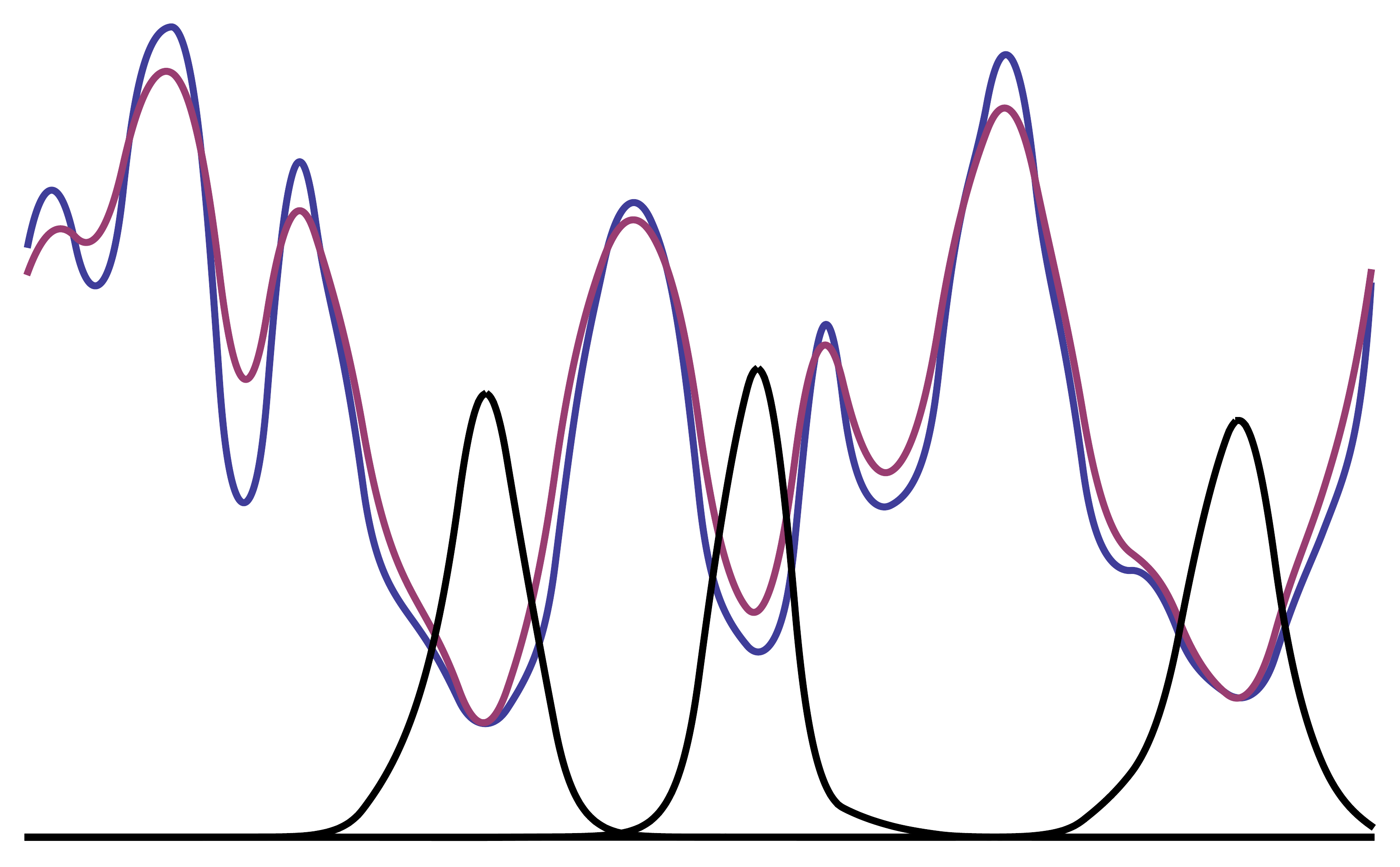}};
\end{tikzpicture}
\caption{The inverse landscape function $1/u$ (blue) and our localized landscape $V_t = V* k_t$ (purple). The first three eigenfunctions (black) localizing in the regularized potential wells.}
\end{figure}
\end{center}

\subsection{Open Problems.} We have proposed a somewhat axiomatic definition of a regularized potential: a local averaging operator that averages the potential $V$ over a fixed scales while affecting the eigenvalue equation
as little as possible. Many questions remain; some of the more obvious questions are as follows:

\begin{enumerate}
\item \textit{Other Kernels.} While our regularization $V_t = V*k_t$ is natural in light of Theorem 1 and Theorem 2, it is an interesting question whether there are other convolution kernels with possibly good properties. Our convolution kernel arises naturally as the first-order term in a more complicated expression; it would be interesting to understand whether higher order terms can be useful in practice (see \S 4.3).
\item \textit{The Choice of Scale.} Having a one-parameter family of regularized kernels $V_t$ has some advantages (for example in allowing us to state Theorem 1 and 2) but, in applications, requires us to set a scale. The proof requires us to have $t \lesssim \|V\|_{L^{\infty}}^{-1}$ but this is clearly not quite necessary. A canonical choice seems to choose $t$ to be roughly comparable to the inverse of the locally averaged potential (which is also the regime to which the proof can be extended without further ado); a better understanding of this parameter would be desirable.
\item \textit{Refined Information.} It is known that the inverse of the landscape function $1/u$ provides an accurate estimate not only for the location of a localized eigenvector but also more accurate estimates for eigenvalues and a more accurate eigenvalue count in terms of Weyl's law. It would be interesting to understand how this compares with $V_t$ (or other suitable regularizations). One might expect that many natural types of regularization (say, convolution with a Gaussian $V * g_t$) would lead to improved local Weyl laws as well -- in some sense, $V$ is a remarkably bad predictor for low-frequency eigenfunctions
in the setting considered here and should be outperformed by $V*g_t$ for many radial, localized kernels $g_t$ that approximate the identity as $t \rightarrow 0$. It is less clear which method (say, taking $1/u$ compared to $V_t$ compared to $V*g_t$) performs best in which setting. 
It would be interesting to have a comprehensive numerical comparison of such methods.
\end{enumerate}

\vspace{-5pt}

\subsection{$1/u$, $V_t$ and general landscapes.} It is clear from the various examples that $1/u$ and $V_t$ are closely related. It is a priori conceivable that they are
both trying to measure the same underlying effective potential through two very different philosophies and happen to agree because they are both somewhat successful. However, Theorem 2 makes this seem unlikely. It is thus an interesting question to try and understand in what way
$1/u$ and $V_t$ are connected. Moreover, they are structurally quite different, one arising as the
solution of a partial differential equation, the other via convolution with a fixed kernel (this could be useful in practice since convolutions with a fixed kernel can be computed rather quickly).
One of the possible interpretations of Theorem 1 and Theorem 2 is that both $1/u$ and $V_t$ do indeed approximate a common underlying structure: this underlying structure, at least
in the proof of Theorem 1 and Theorem 2, is the short-time asymptotic expansion of the associated parabolic equation acting as if the potential was given by $V_t = V*k_t$. \\

\textit{General Landscapes.} If that was indeed the case, then it would be possible to derive effective potentials from equations that look like the landscape function but have different right-hand sides.
We recall that the logic of the landscape function, as used in \cite{arnold, arnold2, fil} and several other papers, is that
$$  \mbox{solve}\quad (-\Delta + V)u = 1 \implies \quad  \mbox{use} \quad \frac{1}{u} \quad \mbox{as effective potential.}$$
The natural generalization would then be as follows: let $f :\Omega \rightarrow \mathbb{R}_{>0}$ be any positive function. We can then consider
the equation $(-\Delta + V)v = f$. Clearly, if $f$ is a constant, we recover
the classical landscape function. If $f$ is any other function, then the assumption of $1/u$ acting as an approximation of $V*k_t$ would predict
$$\boxed{  \mbox{solve}\quad (-\Delta + V)v = f \implies \quad \mbox{use} \quad \frac{f*k_t}{v}  \quad \mbox{as effective potential.}}$$
This is being formally derived in \S 4.2. It may also be useful in practical applications: it is quite conceivable that sometimes one does not have control over the right-hand side $f$ of the equation. As for the value of $t$, we refer to \S 2.5.2.
Clearly, this recovers the classical effective potential coming from the landscape function if $f$ is a constant function. Moreover, if $f$ is a 
slowly varying function, then we already know from abundant experiments that $v$ is really quite localized and would treat $f$ as essentially constant on
small scales thereby effectively reducing the argument to the previous argument: we would expect the approximation to be reasonable.
It thus remains to understand the case when $f$ is itself very rapidly oscillating. \\

\begin{figure}[h!]
\begin{minipage}[l]{.49\textwidth}
\centering
\includegraphics[width = 6cm]{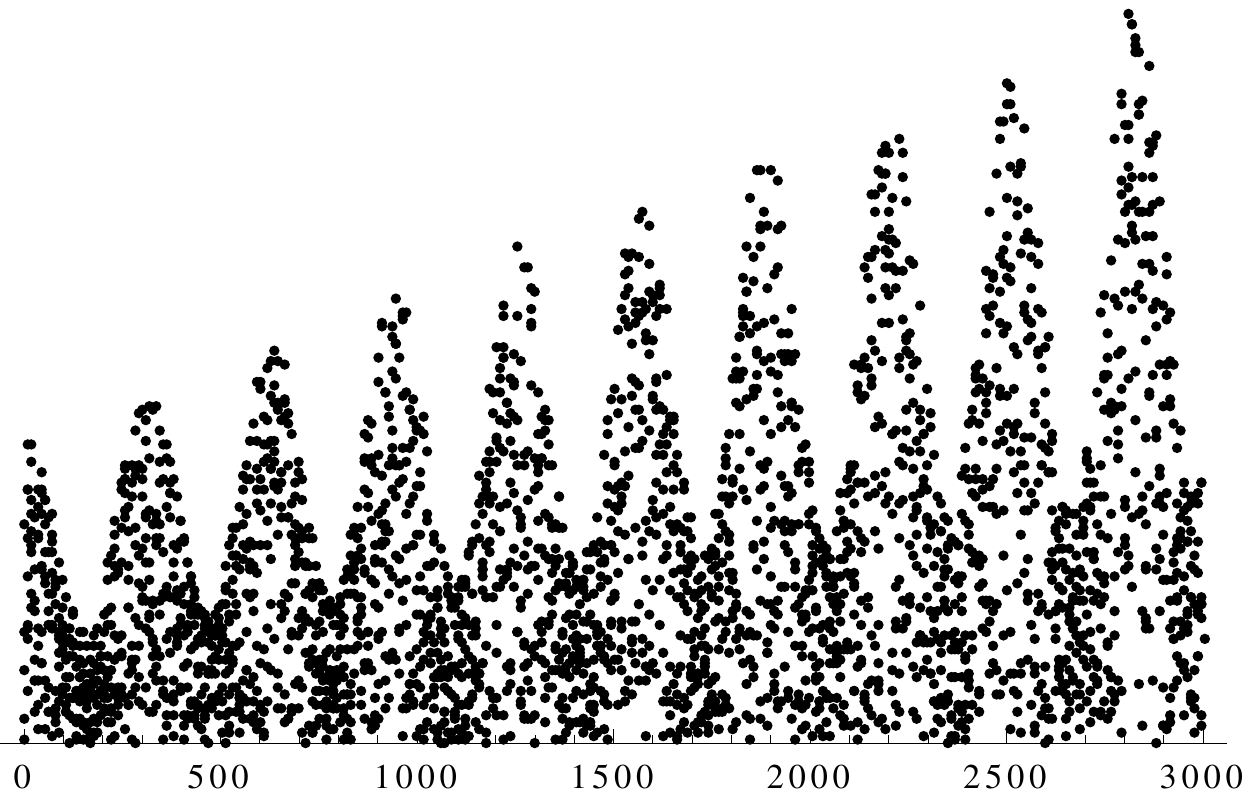} 
\end{minipage} 
\begin{minipage}[r]{.49\textwidth}
\includegraphics[width = 6cm]{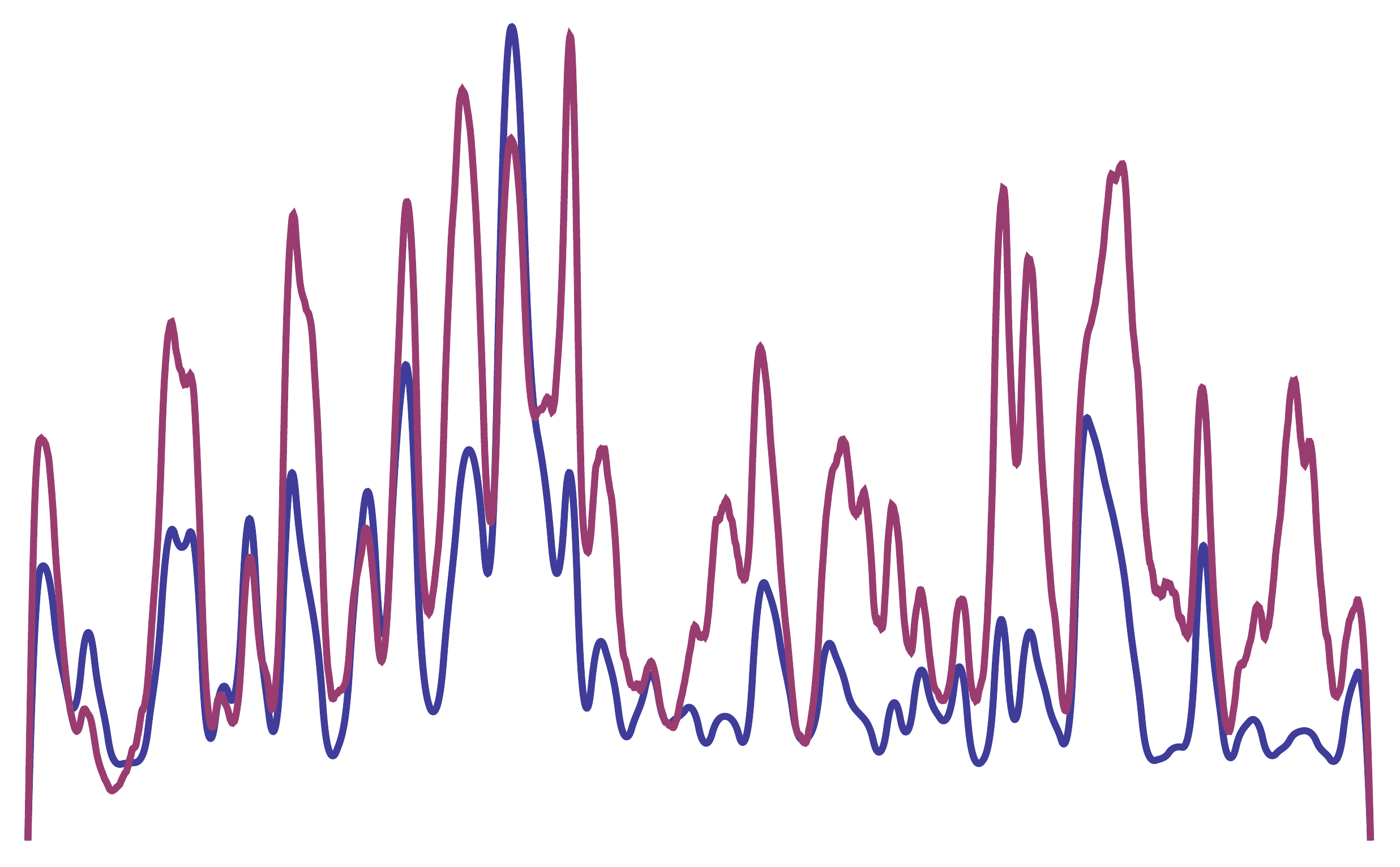} 
\end{minipage} 
\caption{Left: the random vector $f$. Right: the solutions $u$ (blue) and $v$ (purple). $u$ is good at predicting localization, $v$ by itself is not (but see Figure 7 below for $v/(k_t*f)$).}
\end{figure}

Here we considered some simple examples that seem to confirm these ideas.
We choose again the unit interval $[0,1]$ on 3000 nodes with the potential being large (size $\sim 10^5$) and constant on intervals of length $0.01$. We solve for the landscape function
$(-\Delta + V)u=1$ and we solve the equation $(-\Delta + V)v = f$ where $f \in \mathbb{R}^{3000}_{}$ is a randomly generated vector (see Fig. 7) given by $(f_k)_{k=1}^{3000}$ where
$$ f_k = \left(1+\frac{k}{2000}\right)\left(2 + \frac{\cos{(k)}}{50}\right)\cdot \left(\mbox{uniform random variable from}~[0,1]\right).$$
The vector $f$ is shown in the figure above and so are the solutions $u$ and $v$. It is clearly visible that the landscape function $u$ (which we know to be good at
predicting localization) looks very different from $v$ and there is no reason to assume $v$ would have any predictive power. However, once we compute the remaining
correction factor $f*k_t$ and normalize by that function, we recover something strikingly similar to the original landscape function.

\begin{figure}[h!]
\begin{minipage}[l]{.49\textwidth}
\centering
\includegraphics[width = 6cm]{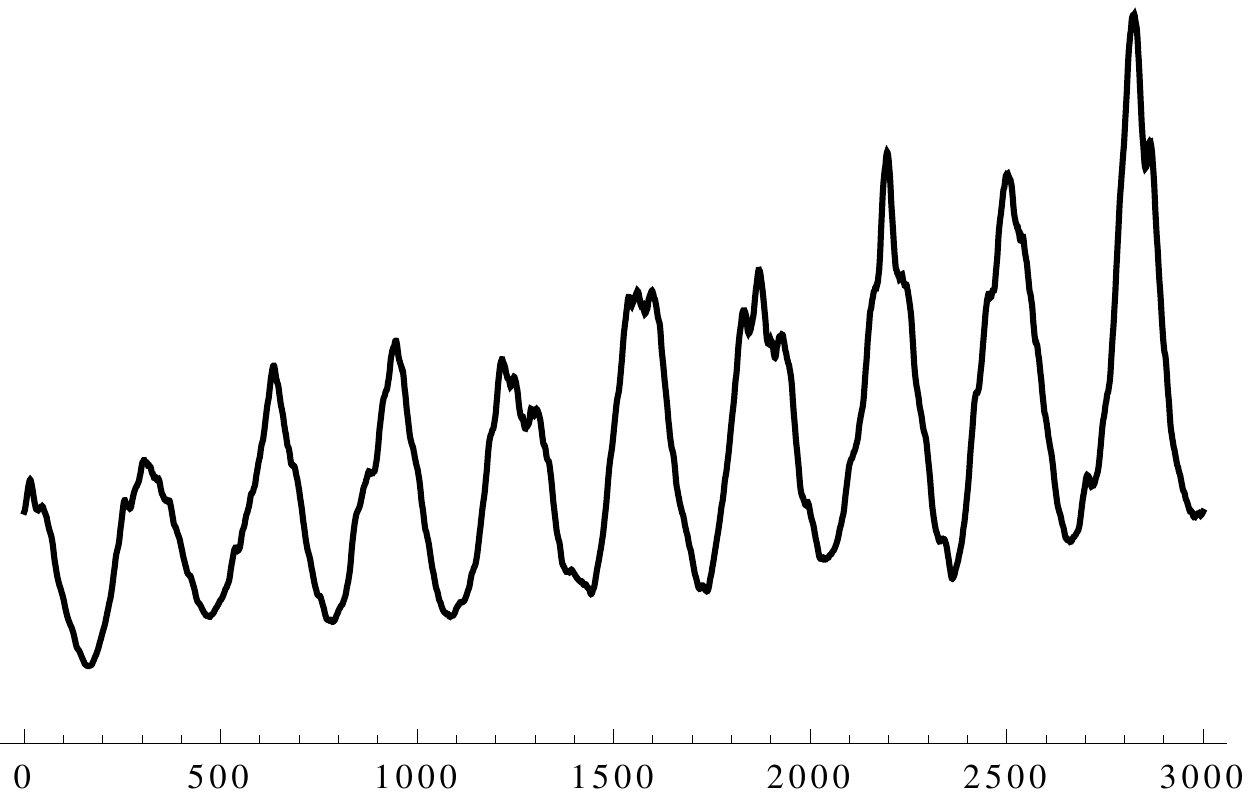} 
\end{minipage} 
\begin{minipage}[r]{.49\textwidth}
\includegraphics[width =6cm]{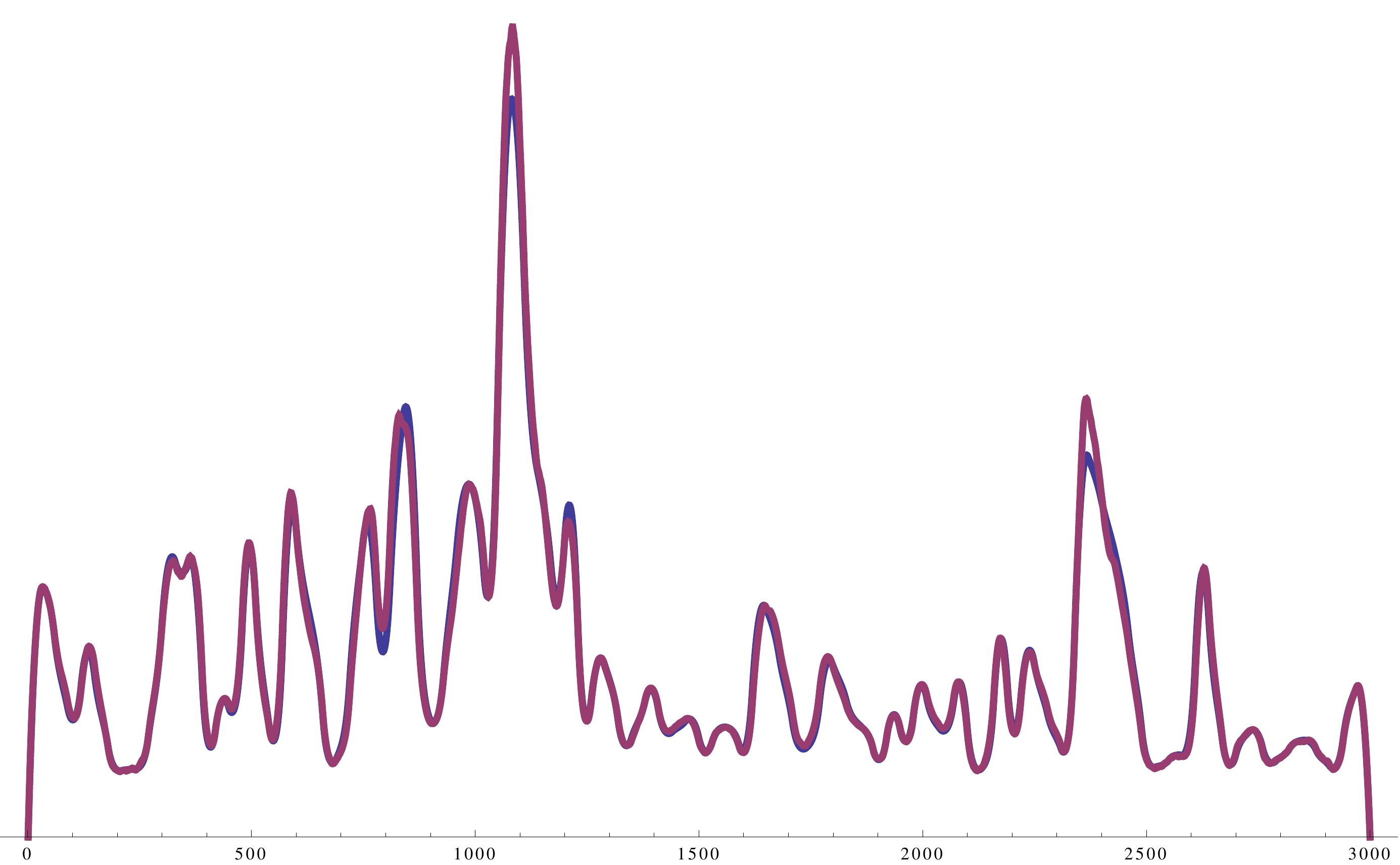} 
\end{minipage} 
\caption{Left: the convolution $f*k_t$ for $t=0.001$. Right: the landscape function $u$ (blue) and $v/(k_t*f)$ (purple).}
\end{figure}

This observed phenomenon seems remarably stable both for random as well as deterministic $f$ and seems to generalize the classical landscape function
to general right-hand sides. Moreover, the generalized landscape with an arbitrary right-hand side $f$ may be useful in practice if $f$ happens to be pre-determined.
This raises a very interesting question when one can choose $f$: is there a particularly clever choice of $f$ that improves on $1/u$ (i.e. the choice $f \equiv 1$) in practice?

\section{Proofs}

\subsection{Ingredients.}
We quickly summarize the main ingredient of our argument for the convenience of the reader. We will associate to the elliptic equation
\begin{align*}
(-\Delta + V)\phi &= \lambda \phi~ \quad \mbox{in~}\Omega \\
 \phi&= 0 \qquad \mbox{on}~ \partial \Omega
\end{align*}
a parabolic equation that we can solve in closed form
\begin{align*}
 \frac{\partial w}{\partial t} &= \Delta w - V w \quad \mbox{in~}\Omega \\
 \omega(0, \cdot)&= \phi(\cdot) \qquad \qquad \mbox{on}~ \partial \Omega.
\end{align*}
Since $\phi$ is an eigenfunction, the equation turns into an ordinary differential equation having the explicit solution
$$ \omega(t,x) = e^{-\lambda t} \phi(x).$$
This, in turn, we will combine with the Feynman-Kac formula which we quickly introduce following the exposition in Taylor \cite{taylor}. We also refer to the book of Lorinczi, Hiroshima and Betz \cite{feyn}.
Following Taylor \cite{taylor}, we have the Feynman-Kac formula in the form
$$ \omega(t,x) = \mathbb{E} \left( \phi(\omega(t)) \exp\left( - \int_{0}^{t} V(\omega(s)) ds\right) \right).$$
Here, the expectation is with respect to Brownian motion $\omega(\cdot)$ started in $x$, running up to time $t$ and being absorbed upon impact on the boundary $\partial \Omega$. We note that we will only
work on $x \in \Omega$ having a positive distance to the boundary $d(x, \partial \Omega) > 0$. Since we will only consider the $t \rightarrow 0^+$ regime, the boundary will not play an important role.
Combining both solutions,
we obtain a reproducing formula
$$ \phi(x) = e^{\lambda t} \cdot \mathbb{E} \left( \phi(\omega(t)) \exp\left( - \int_{0}^{t} V(\omega(s)) ds\right) \right).$$
This formula has already proven useful in a variety of settings \cite{lierl, rachh, stein, steini, steini2}.
We note a special
case of the formula arising from the case $V \equiv 0$: 
$$ \mathbb{E} \left( \phi(\omega(t))  \right) = \mbox{solution of the heat equation at time}~t.$$
In particular, for $t \rightarrow 0$ we have
$$ \mathbb{E} \left( \phi(\omega(t))  \right)  = \phi(x) + t \Delta \phi(x) + \mathcal{O}_{\phi}(t^2),$$
where the implicit constants depends on higher derivatives of $\phi$ in the point $x$. Naturally, this identity could also be derived 
from the explicit form of the heat kernel on domains in $\mathbb{R}^n$: as long as $x$ is sufficiently far away from the boundary,
for short time ($t \ll d(x, \partial \Omega)^2$) and, in particular, as $t \rightarrow 0$, the heat kernel is merely a Gaussian with an exponentially
small error term. The same asymptotic expansion then simply follows from a Taylor expansion of $\phi$ in $x$.
Moreover, a Brownian particle $\omega(t)$ is, in expectation, roughly distance $\sim \sqrt{t}$ away from $x$ and unlikely to travel any further than that.

\subsection{Proof of Theorem 1}

\begin{proof}[Proof of Theorem 1] We use the reproducing formula
\begin{equation} \label{eins}
 \phi(x) = e^{\lambda t} \cdot \mathbb{E}\left( \phi(\omega(t)) \exp\left( - \int_0^t V(\omega(s)) ds\right)\right),
 \end{equation}
 where the expectation is being taken over all Brownian particle starting in $x$, running up to time $t$ and being absorbed at the boundary.
 We will obtain a suitable expansion for $0 \leq t \ll 1$. Clearly,
\begin{equation} \label{zwei}
  e^{\lambda t}  = 1 + \lambda t + \mathcal{O}(t^2).
  \end{equation}
We now compute a first order expansion of the expectation in \eqref{eins} in $t$.
 A Brownian
motion travels distance $\sim \sqrt{t}$ after $t$ units of time and thus, as discussed above,
\begin{equation} \label{drei}
 \mathbb{E} \left(  \phi(\omega(t)) - \phi(x) \right) = t \Delta \phi(x) + \mathcal{O}_{\phi}(t^2)
 \end{equation}
is a first order term. We split the sum into
\begin{align} \label{vier}
\mathbb{E}\left( \phi(\omega(t)) \exp\left( -\int_0^t V(\omega(s)) ds\right)\right) &=  \phi(x)  \mathbb{E}\left(\exp\left(- \int_0^t V(\omega(s)) ds\right)\right) \\
&+ \mathbb{E}\left( ( \phi(\omega(t)) - \phi(x)) \exp\left(- \int_0^t V(\omega(s)) ds\right)\right). \nonumber
\end{align}
A simple Taylor expansion shows for $t$ sufficiently small (depending only on $\|V\|_{L^{\infty}}$),
\begin{equation} \label{funf}
 \mathbb{E}\left(\exp\left(- \int_0^t V(\omega(s)) ds\right)\right)  =  1- \mathbb{E}  \int_0^t V(\omega(s)) ds + \mathcal{O}_{\|V\|_{L^{\infty}}}(t^2).
 \end{equation}
The second term in \eqref{vier} is straightforward to analyze since the first quantity is already of size $\sim t$ and thus, using \eqref{drei},
\begin{align} \label{sechs}
 \mathbb{E}\left( ( \phi(\omega(t)) - \phi(x)) \exp\left(- \int_0^t V(\omega(s)) ds\right)\right) &= t \Delta \phi(x) + \mathcal{O}_{\phi,\|V\|_{L^{\infty}}}(t^2).
\end{align}
Combining \eqref{vier}, \eqref{funf} and \eqref{sechs}, we obtain
\begin{align} \label{sieben}
 \mathbb{E}\left( \phi(\omega(t)) \exp\left( - \int_0^t V(\omega(s)) ds\right)\right) &= \phi(x) - \phi(x) \mathbb{E}  \int_0^t V(\omega(s)) ds \\
 &+t \Delta \phi(x) + \mathcal{O}_{\phi,\|V\|_{L^{\infty}}}(t^2). \nonumber
\end{align} 
Collecting \eqref{eins}, \eqref{zwei} and \eqref{sieben}, we obtain
\begin{align} \label{acht}
 \phi(x) &= \phi(x)  - \phi(x) \mathbb{E}  \int_0^t V(\omega(s)) ds + t \Delta \phi(x) + \lambda t \phi(x) + \mathcal{O}_{\phi,\|V\|_{L^{\infty}}}(t^2).
 \end{align}
Since the left-hand side of \eqref{acht} does not depend on $t$, the coefficient in front of the linear term in \eqref{acht} has to vanish and
\begin{align} \label{neun}
  - \phi(x) \mathbb{E}  \frac{1}{t} \int_0^t V(\omega(s)) ds + \Delta \phi(x) + \lambda \phi(x) = \mathcal{O}_{\phi,\|V\|_{L^{\infty}}}(t).
\end{align}
which we rewrite as
\begin{align} \label{zehn}
 - \Delta \phi(x) +  \phi(x) \mathbb{E}  \frac{1}{t} \int_0^t V(\omega(s)) ds =  \lambda \phi(x) + \mathcal{O}_{\phi,\|V\|_{L^{\infty}}}(t).
\end{align}
 It remains to compute the expectation. In the bulk of the domain, the Brownian particle moves in an unrestricted fashion (as long as $0 \leq t \ll d(x,\partial \Omega)$ and certainly as $t \rightarrow 0$) and we can exchange the order of integration to obtain
\begin{align} \label{elf}
 \mathbb{E} \frac{1}{t} \int_0^t V(\omega(s)) ds  = \frac{1}{t} \int_0^t   \mathbb{E}  V(\omega(s)) ds  
 \end{align}
 The particle $\omega(s)$ in free space is distributed like
 a Gaussian centered at $x$ with variance $s$. Therefore, assuming $t$ is sufficiently small compared to $d(x, \partial \Omega)$,
\begin{align} \label{zwolf}
 \mathbb{E}  V(\omega(s)) = \int_{\Omega}{ V(y)  \frac{ \exp\left( - \|x-y\|^2/ (4s) \right)}{(4 \pi s)^{d/2}}  dy} + \mbox{l.o.t},
 \end{align}
 where the lower order terms are exponentially small provided $t \ll d(x, \partial \Omega)^2$ and converge exponentially in $t$ to 0
 as $t \rightarrow 0$.
 Integrating this in time, we obtain
\begin{align} \label{dreizehn}
\frac{1}{t} \int_0^t   \mathbb{E}  V(\omega(s)) ds = \int_{\Omega}{ V(y)  \frac{1}{t} \int_0^t \frac{ \exp\left( - \|x-y\|^2/ (4t) \right)}{(4 \pi t )^{d/2}} ds  dy} + \mbox{l.o.t.}
\end{align}
 which is the desired statement since the lower order terms in \eqref{zwolf} and \eqref{dreizehn} vanish exponentially when $t \rightarrow 0$.
\end{proof}

\subsection{Proof of Theorem 2}

\begin{proof}[Proof of Theorem 2]
We will, following an idea from \cite{steini}, evaluate the short time evolution of the heat equation
\begin{align} \label{vierzehn}
 \frac{\partial f}{\partial t} &= \Delta f - V f\\
 f(0,x) &= u(x) \nonumber
 \end{align}
 subject to Dirichlet boundary conditions. \eqref{vierzehn} can be solved using the Feynman-Kac formula
 \begin{align} \label{funzehn}
 f(t,x) = \mathbb{E}\left( u(\omega(t)) \exp\left( - \int_0^t V(\omega(s)) ds\right)\right).
 \end{align}
 Arguing as above, we obtain, for $t$ sufficiently small, that
 \begin{align} \label{sechszehn}
  f(t,x) &= u(x) + t \Delta u(x) - u(x) \mathbb{E} \int_0^t V(\omega(s)) ds + \mathcal{O}_{u,\|V\|_{L^{\infty}}}(t^2)\\
  &= u(x) + t \Delta u(x) - u(x) t (V*k_t)(x) + \mathcal{O}_{u,\|V\|_{L^{\infty}}}(t^2)\nonumber
 \end{align}
 However, recalling that
$$ \Delta u - V u = -1$$
we can identify the linear term as
\begin{align} \label{siebzehn}
 f(t,x) = u(x) - t + \mathcal{O}_{u, \|V\|_{L^{\infty}}}(t^2).
 \end{align}
 Comparing \eqref{sechszehn} and \eqref{siebzehn}, we see that
 $$ (-\Delta + V*k_t) u = 1 +  \mathcal{O}_{u, \|V\|_{L^{\infty}}}(t).$$
\end{proof}

\section{Remarks}
\subsection{Unbounded Potentials}
We quickly comment on the condition $V \in C^{}(\overline{\Omega})$ which is not necessarily.  Figuratively put, if $V$ has mild singularities, then these do not affect any of the
arguments as long as the path integral does not become too singular. This is nicely encapsulated in Khasminskii's lemma.

\begin{lemma}[Khasminskii's lemma] Let $V \geq 0$ be a measurable function and $(X_s)_{s\geq 0}$ be a Markov process on $\mathbb{R}^{n}$ with the property that for some $t>0$ and $\alpha < 1$,
$$ \sup_{x \in \mathbb{R}^{d}}{ \mathbb{E}_{x} \left[ \int_0^t V(\omega(s)) ds \right] }= \alpha.$$
Then
$$ \sup_{x \in \mathbb{R}^d}{ \mathbb{E}_{x} \left[ \exp \left( \int_{0}^{t}{V(\omega(s))ds} \right) \right]} \leq \frac{1}{1-\alpha}.$$
\end{lemma}

We refer to the textbook  by Lorinczi, Hiroshima \& Betz \cite{feyn} and the survey article by B. Simon \cite{sim}.  Since we are arguing in the $t \rightarrow 0$ limit, we can essentially relax the condition 
to $V \in L^{n/2, 1}(\Omega)$, where $L^{n/2, 1}$ denotes the Lorentz space (see \cite{lierl} for the details). We especially refer to the nice book \cite{feyn}.

\subsection{Deriving More Landscape Functions.}
The purpose of this section is to explain why the solution of
$$ (-\Delta + V) v = f$$
has a chance of giving rise to an approximation of the effective potential which is approximately $(f*k_t)/v$. We first return to the original landscape function $u(x)$ and write it, as before, as the fixed point of the heat equation
$$ \frac{\partial u}{\partial t} = \Delta v - V v + f.$$
By virtue of the definition of $f$, this equation is stationary in time (much in the same spirit as when we obtained the reproducing identity for Theorem 1). The Feynman-Kac formula allows us to solve this for short time and obtain a reproducing identity
$$ v(x) = v(x) + t\Delta v + t (f*k_t) - t(V*k_t)v + \mbox{l.o.t.}$$
As before, the coefficient in front of the linear term has to vanish leading to
$$ -\Delta v + V_t v = f*k_t + \mbox{l.o.t.}$$
We can now derive the classical asymptotic for the landscape function $u$ by setting $f = 1$, obtaining $f*k_t = 1$ and then assuming that the Laplacian plays a less important role to obtain
$$ V_t \sim \frac{1}{u}.$$
In exactly the same manner, we obtain the more general heuristic
$$ V_t \sim \frac{f*k_t}{v}.$$
As discussed in \S 2.6, this is remarkably well matched by numerical examples.

\subsection{A Second Order Expansion.}
A crucial ingredient of all the arguments in this paper is expansion
\begin{align*}
 \mathbb{E}_{}  \exp \left( \int_0^t V(\omega(s)) ds \right) &= 1 + \mathbb{E}  \int_0^t V(\omega(s)) ds + \mathcal{O}_{\|V\|_{L^{\infty}}}(t^2)\\
 &=1 + t (V*k_t) + \mathcal{O}_{\|V\|_{L^{\infty}}}(t^2)
 \end{align*}
leading to the definition of our kernel $k_t$. However, as is customary in these matters, one could naturally be curious
about the second term in the expansion. We can write
\begin{align*}
\mathbb{E} \left(  \int_0^t V(\omega(s)) ds \right)^2 &= \mathbb{E} \int_{0}^{t} \int_0^t V( \omega(s_1)) V(\omega(s_2)) ds_1 ds_2\\
&= 2\cdot \mathbb{E} \int_{0}^{t} \int_{s_1}^t V(\omega(s_1)) V(\omega(s_2)) ds_2 ds_1\\
&= 2\cdot \int_{0}^{t}\mathbb{E}  \left( V(\omega(s_1))  \int_{s_1}^t  V(\omega(s_2)) ds_2\right) ds_1.
\end{align*}
For any particular Brownian path $\omega(t)$, the behavior of $\omega(s_2)$ for $s_1 \leq s_2 \leq t$ cannot be distinguished from that
of an independent Brownian path started in $\omega(s_1)$ and run for time $t-s_1$. Thus the expectation can be expressed as
\begin{align*}
 \int_{\Omega}  \left( \mathbb{E}  V(\omega(s_1))  \int_{s_1}^t  V(\omega(s_2)) ds_2 \big| \omega(s_1) =y\right) \mathbb{P}(\omega(s)=y)  dx.
\end{align*}
This expectation simplifies and can be written
\begin{align*}
 \left( \mathbb{E}  V(\omega(s_1))  \int_{s_1}^t  V(\omega(s_2)) ds_2 \big| \omega(s_1) =y\right) = V(y) \mathbb{E}  \left( \int_{s_1}^t  V(\omega(s_2)) ds_2 \big| \omega(s_1) =y\right).
\end{align*}
However, we have already evaluated this term in the proof of Theorem 1 and have that 
$$ \mathbb{E}  \left( \int_{s_1}^t  V(\omega(s_2)) ds_2 \big| \omega(s_1) =y\right) = (t-s_1) (V*k_{t-s_1})(y) + \mbox{l.o.t.}$$
Altogether, up to first order, we have
\begin{align*}
\mathbb{E} \left(  \int_0^t V(\omega(s)) ds \right)^2 &= 2 \int_{0}^{t} \int_{\Omega} V(y) \frac{e^{-\|x-y\|^2/(4s)}}{(4 \pi s)^{d/2}} (t-s) (V*k_{t-s})(y) dy ds
\end{align*}
This quantity, an interplay between between Gaussian convolution and convolution with $k_{\cdot}$ then governs the second order
expansion. It remains to be seen whether the second order term can be useful in theory or applications.

\end{document}